\input amstex 
\documentstyle{amsppt}
\input bull-ppt
\keyedby{bull400e/mhm}

\define\reals{\operatorname{\Bbb R}}
\define\e{\varepsilon}
\define\infp{\operatornamewithlimits{inf\vphantom{p}}}
\define\il{\left<}
\define\ir{\right>}
\define\lmp{\text{LMP}}
\define\lmpd{\text{LMP}_d}

\define\ru{R_{\mu}}
\define\mdf{\mu\,(df)}
\define\lall{\Lambda^{\text{all}}}
\define\lstd{\Lambda^{\text{std}}}
\define\Comp{\operatorname{comp}}
\define\Avg{\operatorname{avg}}
\define\Cost{\operatorname{cost}}
\define\co{\Comp^{\Avg}}
\define\er{e^{\Avg}}
\define\a{\alpha}
\define\cost{\Cost^{\Avg}}
\define\app{\text{APP}}

\define\rmin{r_{\min}}
\define\rkhs{\text{reproducing kernel Hilbert space}}
\define\erw{e^{\text{wor}}}
\define\esssup{\operatornamewithlimits{ess\,sup}}
\define\appoo{\text{APP}^{\text{wor}}}

\topmatter
\cvol{29}
\cvolyear{1993}
\cmonth{July}
\cyear{1993}
\cvolno{1}
\cpgs{70-76}
\title Average case complexity\\ 
of linear multivariate problems \endtitle
\author H. Wo\'zniakowski \endauthor
\address Department of Computer Science, Columbia University
and
Institute of Applied 
Mathematics, University of Warsaw\endaddress
\ml henryk\@ground.cs.columbia.edu \endml
\date January 6, 1992\enddate
\thanks This research was
supported in part by the National Science Foundation under 
Contract 
IRI-89-07215 and by AFOSR-91-0347\endthanks
\abstract We study the average case complexity of a linear
multivariate problem $(\lmp)$ defined on functions of $d$ 
variables.
We consider two classes of information. The first $\lstd$ 
consists
of function values and the second $\lall$ of all 
continuous linear
functionals. Tractability of $\lmp$ means that the average 
case
complexity is $O((1/\e)^p)$ with $p$ independent of $d$. 
We prove that
tractability of an $\lmp$ in $\lstd$ is equivalent to 
tractability in
$\lall$, although the proof is {\it not} constructive. We 
provide a
simple condition to check tractability in $\lall$.

We also address the optimal design problem for an $\lmp$ 
by using a relation to the worst case setting. We find the 
order of
the average case complexity and optimal sample points for 
multivariate function approximation. The theoretical 
results are 
illustrated for the folded Wiener sheet 
measure.\endabstract 
\endtopmatter

\document

\heading 1. Introduction \endheading
A linear multivariate problem $(\lmp)$ is defined as 
the approximation of a continuous 
linear operator on functions of $d$ variables. 
Many $\lmp$'s are {\it intractable} in the worst
case setting. That is, the worst case complexity of 
computing 
an $\e$-approximation is infinite or grows {\it 
exponentially} with $d$
(see, e.g., [9]). For example,
consider multivariate integration and function 
approximation of $r$ times
continuously differentiable functions of $d$ variables. 
Then the  worst case
complexity is of order $(1/\e)^{d/r}$
assuming that an $\e$-approximation is computed using 
function values. Thus, if only continuity of the functions
is assumed, $r=0$, then the worst case complexity is 
infinite. For 
positive $r$, if $d$ is large relative to $r$, 
then the worst case complexity is huge even for modest 
$\e$. In either
case, the problem cannot be solved. 

In this paper we study if tractability can be broken by 
replacing 
the worst case setting by an average case setting with a 
Gaussian
measure on the space of functions.  The 
average case complexity is defined as the minimal average 
cost of computing an
approximation with average error at most $\e$. 
We consider two classes of information. 
The first class $\lstd$ consists of function
values, and the second class $\lall$ consists 
of all continuous linear functionals.

We say an $\lmp$ is tractable if the average
case complexity is 
$O\left((1/\e)^p\right)$ with $p$ independent of $d$. The 
smallest such $p$ is
called the exponent of the problem. Under mild 
assumptions, we prove
that tractability in $\lall$ is equivalent to tractability 
in $\lstd$ and that the difference of the exponents is at
most $2$. The proof of this result is {\it not} 
constructive. 
We provide, however, a simple condition to check 
tractability in $\lall$. 

In particular, this means that multivariate integration is 
tractable 
in $\lstd$ and its exponent is at most $2$. 
This should be contrasted with the
worst case setting where, even for $d=1$, the worst case 
complexity 
in $\lstd$ can be infinite or an arbitrary increasing 
function of $1/\e$ (see [14]). 
Of course, intractability of multivariate integration in the
worst case setting can also be broken by switching to the 
randomized
setting and using the classical Monte Carlo algorithm. 

The optimal design problem of constructing sample points 
which achieve 
(or nearly achieve) the average case complexity of an 
$\lmp$ 
in $\lstd$ is a challenging problem. This problem 
has long been open even for multivariate integration and
function approximation. In what follows, we will use the 
word
``optimal'' modulo a multiplicative constant which may 
depend on $d$
but is independent of $\e$. 
Recently, the optimal design problem has been solved 
for multivariate integration for specific Gaussian measures
(see [15] for 
the classical Wiener sheet measure,  
[5] for the the folded Wiener sheet measure, and [13] for 
the
isotropic Wiener measure). 

In this paper, we show under a mild assumption that 
tractability of 
{function approximation} $(\app)$ implies tractability of
other $\lmp$\<s. Therefore, it is enough to address 
optimal sample 
points for $\app$.
Optimal design for $\app$ is analyzed by exhibiting a 
relation between
average case and worst case errors of linear algorithms 
for $\app$.  
This relation reduces the study of the average case to the 
worst case
for a different class of functions. This different class 
is the unit ball
of a reproducing kernel Hilbert space whose kernel 
is given by the covariance kernel of the average case 
measure. 
Similar relations have been used in many papers for 
approximating 
continuous linear functionals; a thorough overview may  be 
found in [11]. 

We illustrate the theoretical results for the folded 
Wiener sheet measure. 
In this case, an $\lmp$ is tractable 
and has exponent at most $2$.  
For $\app$ the exponents in $\lstd$ and $\lall$ are the 
same.
The exponent in $\lall$ was known (see [4]), whereas 
the exponent in $\lstd$ was known to be at most $6$ (see 
[3]). 
Tractability of $\app$ for the folded Wiener sheet measure 
is
in  sharp contrast to intractability of $\app$ for 
the isotropic Wiener measure; see [13].

Tractability of $\app$ in the average case setting 
is significant, since it is known that the randomized 
setting 
does not help (see [12]). Thus, unlike multivariate 
integration, 
intractability of $\app$ in the worst
case setting cannot be broken by the randomized setting.

$\app$ has been studied in $\lstd$ for $d=1$ in [2, 6]. 
For $d\ge 1$, it was shown in [4] that the number of grid 
points needed 
to guarantee an average error $\e$ depends exponentially 
on $d$. 
Of course, $O(\e^{-2-\delta})$ sample points are enough to 
compute an
$\e$-approximation, $\delta>0$. 
Hence, grid points are a poor choice of sample points. 

In [4], the average case complexity of $\app$ in $\lall$ 
was found, 
and it was conjectured that the average case complexity in 
$\lstd$ 
is of the same order.  We prove that this is indeed the 
case.

Optimal design for $\app$ 
is solved by using a relation to the worst
case setting in the reproducing kernel Hilbert space $H$. 
For the folded Wiener sheet measure, $H$ is a
Sobolev space of smooth {\it nonperiodic} functions 
which satisfy certain boundary conditions. 

$\app$ in the worst case setting has been studied
in this Sobolev space
 additionally assuming {\it periodicity} of
functions in [7, 8] (see also [10] for
$d=2)$. It was proven
that {\it hyperbolic cross} points are optimal sample 
points. 
Hyperbolic cross points are defined as a subset of grid 
points
whose indices satisfy a ``hyperbolic'' inequality. 
Approximation
of periodic functions by trigonometric polynomials that 
use Fourier
coefficients with these hyperbolic cross indices was first 
studied in
[1].

For the nonperiodic case, optimal sample
points for $\app$ in the average case
setting are derived from hyperbolic cross points, and the 
average case
complexity is given by
$$
\co(\e;\app)\,=\,\Theta\left(\e^{-1/(\rmin+1/2)}
\big(\log 1/\e\big)^{(k^*-1)(\rmin+1)/(\rmin+1/2)}\right),
$$
with $\rmin=\min_{1\le i\le d}r_i$, 
 where $f^{(r_1,\dots,r_d)}$ is continuous and 
where $k^*$ denotes the number of
$r_i$ equal to $\rmin$. An optimal algorithm 
is given by a linear combination of function values at
sample points derived from hyperbolic cross points. 

Proofs of the results reported here can be found in [16]. 
\heading 2. Linear multivariate problems \endheading
A linear multivariate problem 
$\lmp=\{\lmpd\}$ is a sequence of $\lmpd=(F, \mu$,
$G,S,\Lambda )$
may depend on $d$. We now define them in turn. 

Let $F$ be a separable
Banach space of functions $f:D\to \reals$, $F\subset 
L_2(D)$. 
Here, $D\subset\reals^d$, and its Lebesgue volume
$\l(D)$ is in $(0,+\infty)$. 
We assume that all $L(f)=f(x)$ are in $F^*$.

The space $F$ is equipped with a zero mean Gaussian 
measure $\mu$. 
Let $\ru$ be the covariance kernel of $\mu$, i.e.,
$\ru(t,x)\,=\,\int_Ff(t)\,f(x)\,\mdf\ \text{for } 
t,x\,\in\,D$. 

Let $S:F\to G$ be a continuous linear operator, 
where $G$ is a separable Hilbert space. 
Then $\nu=\mu S^{-1}$ is a zero mean Gaussian measure on 
the Hilbert
space $G$. Its covariance operator $C_\nu=C_\nu^*\ge0$ and 
has a
finite trace.

Finally, $\Lambda$ is either $\lall=F^*$ or $\lstd$ which 
consists of 
$L(f)=f(x),\forall f\in F$, for $x\in D$.

Our aim is to approximate elements $S(f)$ by $U(f)$. 
The latter is defined as follows.
Information about $f$ is gathered by computing a number of 
$L(f)$, where
$L\in\Lambda$, 
$$
N(f)=[L_1(f),L_2(f),\dots,L_n(f)],\quad \forall f\in F.
$$
The choice of $L_i$ and $n=n(f)$
may depend adaptively on the already computed information 
(see [9, Chapter 3]). 
Knowing $y=N(f)$, we compute $U(f)=\phi(y)$ for some 
$\phi:N(F)\to G$. 
The {\it average error} of $U$ is defined as
$$
\er(U)\,=\,\left(\int_F\Vert 
S(f)-U(f)\Vert^2\,\mdf\right)^{1/2}.
$$
To define the average cost of $U$, assume that each 
evaluation of $L(f)$, $L\in\Lambda$ and $f\in F$, costs 
$c=c(d)>0$.
Assume that we can perform arithmetic operations and 
comparisons on real
numbers as well as addition of two elements from $G$ and 
multiplying an element from $G$ by a scalar; all of them 
with cost 
taken as unity. Usually $c\gg 1$. 

For $U(f)=\phi(N(f))$, 
let $\text{cost}(N,f)$ denote the information cost of 
computing
$y=N(f)$. Clearly, we have $\text{cost}(N,f)\ge c n(f)$. 
Let 
$n_1(f)$ denote the number of operations needed to compute 
$\phi(y)$
given $y=N(f)$. (It may happen that $n_1(f)=+\infty$.) 
The {\it average cost} of $U$ is then given as
$$
\cost(U)\,=\,\int_F(\,\text{cost}(N,f)\,+\,n_1(f)\,)\,\mdf.
$$
The {\it average case complexity} of $\lmpd$ is 
the minimal cost of computing $\e$-approxi\-mations,
$$
\co(\e;\lmpd)\,=\,\infp\{\cost(U):\ U\ \text{such that}\ 
\er(U)\le\e\}.
$$
To stress the dependence on certain parameters in 
$\co(\e;\lmpd)$, 
we will sometimes list only those. 
Obviously,  $\co(\e\!;d\!,\lall)\le\co(\e\!;d\!,\lstd)$\<. 
We show that the average case complexity functions in 
$\lall$ and
$\lstd$ are usually closely related. 
\heading 3. Tractability of linear multivariate problems 
\endheading
An $\lmp=\{\lmpd\}$ is called {\it tractable}
if there exists $p\ge 0$ such that for all $d$
$$
\co(\e;\lmpd)\,=\,O\left(c\,\e^{-p}\right).\tag 3.1
$$
The constant in the big $O$ notation may depend on $d$.
The infimum of the numbers $p$ satisfying (3.1) is called 
the {\it exponent} $p^*=p^*(\lmp)$. 
To stress the role of the class 
$\Lambda$, we say that an $\lmp$ is {\it tractable in} 
$\Lambda$ iff
(3.1) holds for $\Lambda$.

In what follows, by {\it multivariate function 
approximation}
we mean $\app=\lmp$ with the embedding $S(f)=I_d(f)=f\in
G=L_2(D)$, where the norm in $L_2(D)$ is denoted by 
$\Vert\cdot\Vert_d$. 

We assume that for all $d$ there exist $K_i=K_i(d)$, 
$i=1,2$, such that
$$\align
&\Vert S(f)\Vert\,\le\,K_1\,\Vert f\Vert_d, \quad 
\forall\,f\in F, 
\tag A.1\\
&\Vert \ru(\cdot,\cdot)\Vert_{L_{\infty}(D)}\,\le\,K_2.
\tag A.2
\endalign
$$
\proclaim{Theorem 3.1}
Suppose {\rm (A.1)} and {\rm (A.2)} hold. 

{\rm (i)} Tractability of $\lmp$ in $\lstd$ is equivalent 
to tractability 
of $\lmp$ in $\lall$ since 
$$
\Comp^{\Avg}(\e;d,\lall)=  
O(c\,\e^{-p(d)}) \text{ implies}\ 
\Comp^{\Avg}(\e;d,\lstd)=O(c\,\e^{-p(d)-2}). 
$$

{\rm (ii)} Let $\lambda_i(d)$ be the ordered eigenvalues 
of the covariance
operator of $\mu S^{-1}$. 
$\lmp$ is tractable in $\lall$  iff 
there exists a positive number $\a$ such that for all $d$,
$$
\sum_{i=n+1}^{+\infty}\lambda_i(d)\,=\,O(n^{-2\a}),
\qquad \text{as}\ \ n \to +\infty. \tag 3.2
$$
The exponent of $\lmp$ is 
$p^*\,=\,1/\sup\{\a:\,\a\ \text{of}\ (3.2)\}$, and
$p^*=+\infty$ if there is no such $\a$. 

{\rm (iii)} Tractability of $\app$ 
in  $\Lambda$ with exponent $p^*$ implies tractability of 
an $\lmp$ in $\Lambda$
with exponent at most $p^*$ provided $\lmp$  differs from 
$\app$ 
only by the choice of $S$. 
\endproclaim

We stress that the proof of Theorem 3.1 is {\it not} 
constructive.
The exponents in $\lall$ and $\lstd$ may differ by at most 
$2$. 
The constant $2$ is sharp. Indeed, for
the integration problem with
the isotropic Wiener measure, the exponent in $\lstd$ 
is $2$ (see [13]), and, obviously, the exponent in
$\lall$ is zero. 
\heading 4. Relation to worst case \endheading
Due to (iii) of Theorem 3.1, 
it is enough to analyze multivariate function approximation 
$\app=\{\app_d\}$ with $\app_d=\{F,\mu,L_2(D),I_d,\lstd\}$. 
The {\it average} case errors of $\app$ are related to
{\it worst} case errors of the same
$I_d$ restricted to a specific subset of $F$. This 
specific subset of $F$ is 
the unit ball $BH_\mu$ of a $\rkhs$ $H_\mu$. 
The space $H_\mu$ is the completion of finite-dimensional 
spaces of the form
$$
\text{span}\left(\ru(\cdot,x_1),\ru(\cdot,x_2),\dots,\ru(%
\cdot,x_k)\right)\.
$$
The completion is with respect to 
$\Vert\cdot\Vert_\mu=\il \cdot,\cdot\ir_\mu^{1/2}$, where 
$\il R(\cdot,x),R(\cdot,t)\ir_\mu=R(x,t)$. 

Consider a linear $U$ which uses sample points $x_j$. That 
is, we have 
$U(f)\,= %
\sum_{j=1}^nf(x_j)\,g_j$, where $g_j\in L_{\infty}(D)$. 
It is easy to show that
$$
\er(U)\,=\,\er(U;\app_d)\,=\,\left(\int_D\Vert
h^*(\cdot,x)\Vert_\mu^2\,dx\right)^{1/2},
$$
where 
$h^*(\cdot,x)\,=\,\ru(\cdot,x)\,-\,\sum_{j=1}^ng_j(x)\,%
\ru(\cdot,x)\ 
\in H_\mu$. 

Consider now the same $U$ for multivariate function 
approximation
in the $L_{\infty}(D)$ norm 
$$
\appoo_d=\{BH_\mu,L_{\infty}(D),I_d,\lstd\}
$$
in the worst case setting. We now assume that $H_\mu$ is a 
subset of 
$L_{\infty}(D)$ and that the embedding $I_d$ maps $H_\mu$ 
into 
$L_{\infty}(D)$. The worst error of $U$ is equal to
$$
\erw(U;\appoo_d)\,=\,\sup\big\{\big \Vert f- U(f)\big
\Vert_{L_{\infty}(D)}: \ \Vert f\Vert_\mu\le1 \big\}. 
$$
It is easy to show that 
$\erw(U;\appoo_d)\,=\,\esssup_{x\in D}\,\Vert 
h^*(\cdot,x)\Vert_\mu$, 
which yields 
$$
\er(U;\app_d)\,\le\,{\sqrt{l(D)}}\,\erw(U;\appoo_d),\tag 4.1
$$
where $l(D)$ is the Lebesgue volume of $D$. 
\heading 5. Application for folded Wiener sheet measures 
\endheading
We assume that $\mu$ is the folded Wiener sheet 
measure (see [4]). That is, $D=[0,1]^d$ and 
$F$ is the space of $r_i$ times continuously
differentiable functions with respect to $x_i$ which 
vanish with their
derivatives at points with at least one component equal to 
zero. 
The norm of $F$ is the sup norm
on $(r_1,\dots,r_d)$ derivatives. 
The covariance kernel $\ru$ of $\mu$ is  
$$ 
\ru(t,x)\,=\,\prod_{j=1}^{d}\int_0^1 
\frac{(t_j-s)^{r_j}_{+}}{r_j!}\,
\frac{(x_j-s)^{r_j}_{+}}{r_j!}\,ds. 
$$
Observe that $\ru(t,t)\le 1$ and
(A.2) holds with $K_2\le 1$.

The space $H_\mu$ consists now of functions $f$ of the 
form (see [5]) 
$$
f(x)=\int_D\prod_{j=1}^d\frac{(x_j-t_j)^{r_j}_+
}{r_j!}\,\phi(t_1,t_2,
\dots,t_d)\,dt_1\,dt_2\cdots dt_d,\quad \forall x\!\in 
\!D,\ 
\phi \!\in \!L_2(D).  
$$
The inner product of $H_\mu$ is 
$\il f,g \ir_\mu\,=\,\int_Df^{(r_1,\dots,r_d)}(t)\,
g^{(r_1,\dots,r_d)}(t)\,dt$. 

Average case errors for $\app_d$ can be bounded (see 
(4.1)) by analyzing 
the worst case of 
$$
\appoo_d\,=\,\{BH_\mu,L_\infty(D),I_d,\lstd\}. 
$$

Let $W_0$ 
be a subspace of $H_\mu$ of periodic functions 
for which $f^{(i_1,\dots,i_d)}\,(t)=0$ for all $i_j\le 
r_j$ and all $t$ 
from the boundary of $D$. 
Multivariate function approximation for the unit ball 
of $W_0$ in the worst 
case setting has been analyzed by 
Temlyakov in [7, 8].  He constructed sample points $x_j$ 
and functions
$a_j$ such that for 
$T_n(f,x)\,=\,\sum_{j=1}^nf(x_j)\,a_j(x)$ we have
$$
\Vert f-T_n(f,\cdot)\Vert_{L_{\infty}(D)}\,=\,
O(n^{-(\rmin+1/2)}\,(\log n)^{(k^*-1)(\rmin+1)}
),\tag 5.1
$$
where $\rmin\,=\,\min\{r_j:\,1\le j\le d\}\ \text{and}\ 
k^*\,=\,\text{card}(\{j:\ r_j=\rmin\})$. 

The sample points $x_j$ are called {\it hyperbolic cross 
points} 
and the functions $a_j$ are obtained by linear 
combinations of the 
{\it de la Vall\'ee-Poussin kernel}. 

To extend Temlyakov's result to nonperiodic functions, 
define for $f$ from $BH_\mu$ 
$$
g(x)\,=\,f\big(\vec h(x)\big),\quad \forall x\in D,
$$
where $\vec h(x)=\,\big(h(x_1),h(x_2),\dots,h(x_d)\big)$ 
and 
$h(u)=\,4\,u\,(1-u),\ \forall u\in [0,1]$. 

Observe that $g$ is periodic and enjoys the same 
smoothness as $f$; 
that is, $g\in W_0$. 
There exists a constant $K=K(d,\vec r)$ such that
$\Vert g\Vert_\mu\,\le\,K$.
Define
$$
U^*_n(f,t)\,=\,T_n(g,\vec h^{-1}(t)),
$$
where $\vec h^{-1}(t)=(\tfrac12(1-\sqrt{1-t_1}),\dots,
\tfrac12(1-\sqrt{1-t_d}))$, $t\in D$. 
We have 
$$
U^*_n(f,t)\,=\,\sum_{j=1}^nf\big(\vec h(x_j)\big)\,
a_j\big(\vec h^{-1}(t)\big)
\,=\,\sum_{j=1}^nf(x^*_j)h^*_j(t),\tag 5.2
$$
where $x^*_j=\vec h(x_j)$, with a hyperbolic cross point 
$x_j$,
and $h^*_j(t)=a_j\big(\vec h^{-1}(t)\big)$. 

It is possible to check that for all $f$ from  $BH_\mu$ we 
have 
$$
\Vert f\,-\,U^*_n(f,\cdot)\Vert_{L_{\infty}(D)}\,=\,
O(n^{-(\rmin+1/2)}\,(\log n)^{(k^*-1)(\rmin+1)}
).\tag 5.3
$$
From (5.3) and (4.1) we conclude that 
$$
\co(\e;\app_d)\,=\,
O(c\,\e^{-1/(\rmin+1/2)}(\log 
1/\e)^{(k^*-1)(\rmin+1)/(\rmin+1/2)}).\tag 5.4
$$
Clearly, $\co(\e;\app_d)$ is bounded from below by the
corresponding average case complexity in the class 
$\lall$. The latter 
was determined in [4]. These two average case
complexity functions differ by at most a constant. 
Thus, the $O$ in (5.4) can
be replaced by $\Theta$. Furthermore,  
the linear approximation $U^*_n$ given by (5.2) is optimal, 
i.e., $U^*_n$ computes an $\e$-approximation with 
the average cost $(c+2)n$ which is minimal, modulo a 
constant, if 
$$
n\,=\,O(\e^{-1/(\rmin+1/2)}\,(\log
1/\e)^{(k^*-1)(\rmin+1)/(\rmin+1/2)}). \tag 5.5
$$
\proclaim{Theorem  5.1}
For \!$\app$\! the average
case complexity functions $\Comp^{\Avg}(\e\!;d\!,\lstd)$ 
and 
$\Comp^{\Avg}(\e;d,\lall)$ 
differ at most by a constant and 
$$
\Comp^{\Avg}(\e;d,\lstd)\,=\,
\Theta(c\,\e^{-1/(\rmin+1/2)}\,
(\log 1/\e)^{(k^*-1)(\rmin+1)/(\rmin+1/2)}).
$$
The linear $U^*_n$ given by {\rm (5.2)}
which uses $n$ sample points derived from the hyperbolic 
cross points with
$n$ given by {\rm (5.5)} is optimal in the classes $\lstd$ 
and $\lall$. 
\endproclaim

From Theorem 5.1 we have that $\app$
is tractable in $\lstd$ since $1/(\rmin+1/2)\le 2$. The 
exponent of
$\app$ is the same in $\lall$ and $\lstd$. Since $r_i$ may 
depend on
$d$, we have 
$$
p^*(\lstd)\,=\,\left(1/2\,+\,\min\{r_j(d):\ j=1,2,\dots,d\ 
\text{and}\ 
d=1,2,\dots\}\right)^{-1}\ \le\ 2.
$$
Obviously, any $\lmp$ which satisfies (A.1) and which is
equipped with the folded Wiener sheet measure 
is tractable and has exponent at most $p^*(\lstd)\le 2$.
\heading Acknowledgment \endheading
I thank A. Papageorgiou, S. Paskov, L. Plaskota,
V. N. Temlyakov,  J. F. Traub, G. W. Wasilkowski, and A. 
G. Werschulz  for
valuable comments.

\enddocument